\documentclass{article}
 \usepackage{amssymb}
\usepackage{amsmath}
\title{Integral and isocapacitary inequalities}
\author{Vladimir Maz'ya\footnote{The author was partially supported by  the UK Engineering and Physical
Sciences Research Council grant EP/F005563/1.}
\\*[10pt]
\emph{\small Department of Mathematical Sciences,   University of Liverpool, Liverpool L69 7ZL}
\\
\emph{\small and}\\
\emph{\small Department of Mathematics, Link\"oping University, Link\"oping, SE-581 83}
\\
\emph{\small e-mail: vlmaz@mai.liu.se}
\\[15pt]
\emph{To Victor Havin on   his 75th birthday}
\\[5pt]
\emph{  with 45-year affection and admiration}}
\date{}

\begin{document}
\maketitle

\begingroup

\narrower\noindent
{\it Abstract.}  It is shown by a counterexample that isocapacitary and isoperimetric constants of a multi-dimensional Euclidean domain starshaped with respect to a ball are not equivalent. Sharp integral inequalities involving the harmonic capacity which imply Faber-Krahn property of the fundamental Dirichlet-Laplace eigenvalue are obtained. Necessary and sufficient conditions ensuring integral inequalities between a difference seminorm and the $L_p$-norm of the gradient are found.
\endgroup

\bigskip\noindent
{ {\small\it  Mathematics Subject Classification (2000)}: }{\small Primary 54C40, 14E20; Secondary 46E25, 20C20}

\medskip\noindent
{ {\small\it Key words}: }{\small Dirichlet problem, Laplace operator, eigenvalue, isoperimetric  inequalities, isocapacitary inequalities, area minimizing function, Besov spaces, difference seminorms

\section{Introduction}

Isocapacitary inequalities are intimately connected with various properties of Sobolev spaces, especially with norms of embedding operators \cite{Gr}, \cite{M1}--\cite{M4}, \cite{M6}, \cite{M7}, \cite{M9}, \cite{M14}. For instance, the best constants in some of these inequalities give two-sided estimates for eigenvalues of boundary value problems \cite{M3}, \cite{M4},  \cite{M6}, \cite{M7}. Recently, in \cite{CM} and \cite{ACMM},  isocapacitary inequalities were applied to the theory of quasi-linear second order elliptic equations. 

\smallskip

The present paper deals with three topics related to isocapacitary inequalities. First we show by a counterexample in Sect. 2 that the fundamental eigenvalue of the Dirichlet Laplacian is not equivalent to  an  isoperimetric constant, called, as a rule, Cheeger's constant \cite{Ch}, in contrast with an isocapacitary constant introduced in \cite{M3} (see also \cite{M6}).\footnote{ By the equivalence of the set functions $a$ and $b$, defined on subsets of $ \Bbb{R}^n$, I mean the existence of positive constants $c_1$ and $c_2$ depending only on $n$ and such that $c_1 a \leq b\leq c_2 a$.} This equivalence, even uniform with respect to the dimension,  holds for convex domains as proved  recently by B. Klartag and E. Milman (oral communication) but, as I show,  it   fails even in the class of domains starshaped with  respect to a ball.

\smallskip

 Sect. 3 is devoted to certain integral capacitary inequalitites which are stronger than the classical Faber-Krahn property of the fundamental Dirichlet-Laplace eigenvalue (see \cite{PS}).

\smallskip

In Sect. 4 and  5 one can find necessary and sufficient conditions for the inequality
\begin{equation}\label{k00}
\Bigl( \int_\Omega\int_\Omega |u(x) -u(y)|^q \mu(dx, dy)\Bigr)^{1/q} \leq C\, \|\nabla u \|_{L_p(\Omega)}
\end{equation}
formulated in terms of the isoperimetric ($q\geq p=1$) and isocapacitary ($q>p>1$) inequalities of a new type. Here $\Omega$ is a subdomain of a Riemannian manifold, $\mu$ is a given measure of   two subsets of $\Omega$ and $u$ is an arbitrary smooth function.  

\smallskip

No caracterization of (\ref{k00}) was known previously even for functions on the real line $\Bbb{R}$ (see Problem 3 in \cite{KMP}). A particular case of a result obtained at the end of Sect. 5  is the criterion of the validity of (\ref{k00}) for all $u\in C^\infty_0(\Bbb{R})$:
$$
\mu \bigl(([\alpha, \beta], \Bbb{R} \backslash (\alpha -r, \beta +r)\bigr)\leq const.\, r^{-q(1-p)/p},$$
where $r>0$, $\alpha<\beta$, and the constant factor does not depend on $\alpha$, $\beta$, and $r$.

\smallskip

The marginal value  $q=p>1$ in (\ref{k00}) has special features and  a sufficient condition for (\ref{k00}) is given  in Sect. 6.  The article is finished  with a short discussion of the inequality
$$\Bigl( \int_\Omega |u|^q\, \nu(dx) \Bigr)^{1/q} \leq \Bigl( \int_\Omega\int_\Omega |u(x) -u(y)|^p \mu(dx, dy)\Bigr)^{1/p}
$$
with a nonnegative measure $\nu$ in $\Omega$, $\mu$ as above, and $q\geq p\geq 1$.

\smallskip

It is worth mentioning that  a Riemannian structure of $\Omega$ is not very important for most of the results presented  in Sect. 3-5.  It can be   replaced by some natural requirements on the $p$-energy integral on a metric  space (see \cite{M12}, \cite{KK}). 

\smallskip

In this article, I use a number of assertions from the book \cite{M10} which are not formulated in detail but supplied with  references. It is therefore helpful to read the paper with \cite{M10} close at hand.

\section{The first Dirichlet-Laplace eigenvalue\\  and an isoperimetric constant}

Let $\Omega$ be a subdomain of a $n$-dimensional Riemannian manifold $\frak{R}_n$ and 
let $\Lambda(\Omega)$ be the first eigenvalue of the Dirichlet priblem for the Laplace operator $-\Delta$ in $\Omega$ or,
more generally, the upper lower bound of the spectrum of this operator:
\begin{equation}\label{k1}
\Lambda(\Omega) = \inf\limits_{u\in C_0^\infty(\Omega)} \frac{\|\nabla u\|^2_{L_2(\Omega)}}{\|u\|^2_{L_2(\Omega)}}.
\end{equation}
By \cite{M3} (see also \cite{M6} and Corollary 2.3.3 \cite{M10}), $\Lambda(\Omega)$ admits the two-sided estimate
\begin{equation}\label{k2}
\frac{1}{4}\Gamma(\Omega) \leq \Lambda(\Omega)\leq \Gamma(\Omega)
\end{equation}
with
$$\Gamma(\Omega): =\inf\limits_{\{F\}} \frac{{\rm cap}(F; \Omega)}{m_n(F)}.$$
By $m_n$ the $n$-dimensional Lebesgue measure on $\frak{R}_n$ is meant,  the infimum is taken over all compact subsets of
$\Omega$ and ${\rm cap}(F; \Omega)$ stands for the relative harmonic capacity of $F$ with respect to $\Omega$:
$${\rm cap}(F; \Omega) = \inf\Bigl\{\int_\Omega |\nabla u|^2 dx: \, u\in C_0^\infty(\Omega), \, u\geq 1 \,\, {\rm on}\,\, F\Bigr\}.$$
We write ${\rm cap}(F)$ instead of ${\rm cap}(F; \Bbb{R}^n)$.

\smallskip

By Theorem 2.2.1 \cite{M10},  the set function
\begin{equation}\label{k3}
\gamma(\Omega) = \inf\limits_{u\in C_0^\infty(\Omega)} \frac{\|\nabla u\|_{L_1(\Omega)}}{\|u\|_{L_1(\Omega)}} 
\end{equation}
admits the geometric representation
\begin{equation}\label{k4}
\gamma(\Omega) =\inf\limits_{\{g\}} \frac{H_{n-1}(\partial g)}{m_n(g)},
\end{equation}
where $g$ is an arbitrary open subset of $\frak{R}_n$ with compact closure $\overline{g}$  in $\Omega$ and smooth boundary $\partial g$, and $H_{n-1}$ is the
$(n-1)$-dimensional Hausdorff measure. Obviously, for all $u\in C_0^\infty(\Omega)$,
$$\gamma(\Omega) \leq \frac{\displaystyle{\int_\Omega |\nabla (u^2)| dx}}{\displaystyle{\int_\Omega u^2 dx}} \leq 2\,
\frac{\|\nabla u\|_{L_2(\Omega)}}{\|u\|_{L_2(\Omega)}}.$$
Hence
\begin{equation}\label{k5}
\gamma(\Omega)^2 \leq 4\, \Lambda(\Omega), 
\end{equation}
which shows, together with (\ref{k2}) and (\ref{k3}), that
\begin{equation}\label{k6}
\gamma(\Omega)^2 \leq 4\, \Gamma(\Omega)
\end{equation}
(the square of the isoperimetric constant is dominated by the isocapacitary one).

\smallskip

One can ask whether an upper bound for $\Gamma(\Omega)$ formulated in terms of $\gamma(\Omega)$ exists.
The negative answer is obtained easily if the class of sets $\Omega$ is not restricted. In fact, let $F$ be a compact
subset of the open $n$-dimensional unit cube $Q$ in the Euclidean space $\Bbb{R}^n$, such that
$$H_{n-1}(F) =0 \qquad {\rm and} \qquad  {\rm cap} (F) >0.$$
 By $\Omega$ we shall mean the complement of the union of all integral shifts of $F$.
 Now, by Theorem 11.2 \cite{M10}, $\gamma(\Omega) =0$ and $\Gamma(\Omega) \geq \Lambda(\Omega) >0$, which
gives the  negative answer to the question formulated above.

 \smallskip

 Let us put  the same question  for  domains in $\Bbb{R}^n$ starshaped with respect to balls. We   show that the answer
stays negative in a certain sense.

 \medskip

 {\bf Example.} Let $\Omega$ be a subdomain of the $n$-dimensional unit ball $B$, starshaped with respect to a
concentric ball $B(0; \rho) = \{x: |x|<\rho\}$. Here we show that the inequality opposite to (\ref{k6}):
\begin{equation}\label{k7}
\gamma(\Omega)^2 \geq C\, \Gamma(\Omega)
\end{equation}
is imposssible with $C$ independent of $\rho$.  Moreover, we shall construct a sequence of domains
$\{\Omega_N\}_{n\geq 1}$ situated in $B$ and such that

$(i)$ $\Omega_N$ is starshaped with respect to a ball $B(0, \rho_N)$, where $\rho_N\to 0$,

$(ii)$ $\Gamma(\Omega_N) \to \infty$,

$(iii)$ $\gamma(\Omega_N)\leq c$, where $c$ depends only on $n$.

\smallskip

Let $N$ stand for a sufficiently large integer number. By $\{\omega_j\}_{j= 1}^{N^{n-1}}$ we denote a collection of
points on the unit sphere $S^{n-1}$ placed  uniformly   in the sense that the distance from every point $\omega_j$ to  the set of 
other points of the collection lies between $c_1\, N^{-1}$ and $c_2\, N^{-1}$, where $c_1$ and $c_2$ are 
positive constants, depending only on $n$. Consider a closed rotational cone $C_j$ with the axis $O\omega_j$ and the vertex at the distance
$c_0 N^{-1}$ from $O$, where $c_0$ is an absolute constant large enough. The opening of $C_j$ will be independent
of $j$ and denoted by $\varepsilon_N$. Let $\varepsilon_N = o(N^{\frac{1-n}{n-2}})$. 
Clearly, the complement of $C_j$ is visible from a sufficiently small  ball $B(0; \rho_N)$.
Therefore, the domain
$$\Omega_N : = B\backslash  \displaystyle{{\cup}_{j}}
C_j$$
is starshaped with respect to $B(0, \rho_N)$.

\smallskip

We shall find the limit of $ \gamma(\Omega_N)$ as $N\to \infty$  as well as a lower estimate for $\Gamma(\Omega_N)$.
Clearly, $\gamma(\Omega_N) \geq \gamma(B) =n$. Furthermore, by (\ref{k4}),
\begin{eqnarray*}
\gamma(\Omega_N) &\leq& \frac{H_{n-1}(\partial \Omega_N)}{m_n(\Omega_N)} =
\frac{H_{n-1}(\partial B) + H_{n-1}\bigl( \cup_j (B\cap\partial C_j)\bigr)}{m_{n}(B) - m_{n}\bigl( \cup_j (B\cap C_j)\bigr)}\\
&\leq & \frac{|S^{n-1}| + c_1\,  \varepsilon_N^{n-1}\, N^{n-1}}{|S^{n-1}|/n - c_2\, \varepsilon_N^{n-2}\, N^{n-1}}
\end{eqnarray*}
and therefore, 
$$\lim\limits_{N\to \infty} \gamma(\Omega_N) =n.$$

\smallskip

In order to estimate $\Gamma(\Omega_N)$ from below, we construct a covering of $B$ by the balls ${\cal B}_k:= B(x_k, 4c_0N^{-1})$,
whose multiplicity does not exceed a  constant depending only on $n$.  Let $|x_k| \geq c_0\, N^{-1}$. Theorem 10.1.2  \cite{M10} implies
 \begin{equation}\label{k8}
 c\, N^n {\rm cap}({\cal B}_k\backslash \Omega_N)\int_{{\cal B}_k} u^2\,  dx \leq \int_{{\cal B}_k} |\nabla u|^2 dx
 \end{equation}
 for all $u\in C_0^\infty(\Omega_N)$, and the result will stem from a proper lower bound for ${\rm cap} ({\cal B}_k\backslash \Omega_N)$.
 
 \smallskip
 
First, let us consider  $n=3$. Clearly, ${\cal B}_k\backslash \Omega_N$ contains a right rotational  cylinder $T_k$ with height $c_0\, N^{-1}$ and diameter of the base $\varepsilon_N \, N^{-1}$. Now,  by Proposition 9.1.3/1 \cite{M10},
$${\rm cap} (T_k) \geq c\, N^{-1}|\log \varepsilon _N|^{-1}.$$
This estimate in combination with (\ref{k8})  gives

\begin{equation}\label{k9}
c\, N^2 |\log \varepsilon_N|^{-1}\int_{{\cal B}_k} u^2\,  dx \leq  \int_{{\cal B}_k} |\nabla u|^2 dx.
 \end{equation}
 Choosing $\varepsilon_N = {\rm exp} (-N)$ and summing (\ref{k9}) over all balls ${\cal B}_k$, we obtain
$\lambda(\Omega_N)\geq c\, N$. Hence $\lambda(\Omega_N)\to \infty$ where as $\gamma(\Omega_N) \leq c$.
Thus, in particular, there is no inequality
 $$\Bigl( \inf\limits_{\{g\}} \frac{s(\partial g)}{m_3(g)}\Bigr) ^2 \geq C \inf\limits_{\{F\}} \frac{{\rm cap}(F;\Omega)}{m_3(F)}$$
 and, equivalently,
 $$\Bigl( \inf\limits_{\{g\}} \frac{s(\partial g)}{m_3(g)}\Bigr) ^2 \geq C \Lambda(\Omega)$$
 with  constant  factors $C$ independent of the radius $\rho$.

\smallskip

For dimensions greater than $3$, the very end of the argument remains intact but  the estimation  of ${\rm cap}({\cal B}_k\backslash \Omega_N)$ becomes a bit more complicated and the choice of $\varepsilon_N$ will be different.

\smallskip

Let $\alpha\, {\cal B}_k$ stand for the ball concentric with ${\cal B}_k$ and dilated with coefficient $\alpha$.  We introduce the set $s_k = \{j: C_j\cap \frac{1}{2} {\cal B}_k \neq \emptyset\}$.  With every $j$ in $s_k$ we associate a right rotational cylinder $T_j$ coaxial with   the cone $C_j$ and   situated in $C_j\cap \frac{1}{2} {\cal B}_k$. The height of $T_j$ will  be equal to $c_0\, N^{-1}$ and the diameter of the base   equal to $\varepsilon_N\, |x_k|$. We define a cut-off function $\eta_j$, equal to $1$ on the $\varepsilon_N\, |x_k|$-neighbourhood of $T_j$, zero outside the $2\varepsilon_N\, |x_k|$-neighbourhood of $T_j$ and satisfying the estimate
$$|\nabla\eta_j(x)| \leq c\, \delta(x)^{-1},$$
where $\delta(x)$ is the distance from $x$ to the intersection of $ T_j$ with 
the axis of $ C_j$. 

\smallskip

By ${\cal P}_k$ we denote the equilibrium potential of ${\cal B}_k \backslash \Omega_N$. We have
$$\sum_{j\in s_k} {\rm cap}(C_j\cap  {\cal B}_k) \leq  \sum_{j\in s_k} \int_{\Bbb{R}^n}  |\nabla({\cal P}_k\eta_j)|^2 dx$$
$$\leq c\Bigl(\int_{\Bbb{R}^n} |\nabla {\cal P}_k|^2 dx + \int_{\Bbb{R}^n} {\cal P}_k^2 \delta^{-2} dx \Bigr).$$
Changing the constant $c$, one can majorize the last integral by the previous one due to Hardy's inequality. Hence,
\begin{equation}\label{k80}
{\rm cap} ( {\cal B}_k \backslash \Omega_N) \geq c\sum_{j\in s_k} {\rm cap}(T_j).
\end{equation}

\noindent
By Proposition 9.1.3/1 \cite{M10},
$${\rm cap}(T_j) \geq c\, \bigl( \varepsilon_N |x_k|\bigr)^{n-3} N^{-1}.$$
Furthermore, it is visible that the number of integers in $s_k$ is between two multiples of  $|x_k|^{1-n}$. Now, by (\ref{k80})
$${\rm cap} ( {\cal B}_k \backslash \Omega_N) \geq c \, |x_k|^{1-n}\, \bigl( \varepsilon_N |x_k|\bigr)^{n-3} N^{-1}$$
 and by (\ref{k8})
\begin{equation}\label{k81}
c\, N^{n-1} |x_k|^{-2} \varepsilon_N^{n-3} \int_{{\cal B}_k} u^2 dx \leq \int_{{\cal B}_k} |\nabla u|^2 dx.
\end{equation}
Since $|x_k| \leq 1$, it follows by summation of (\ref{k81}) over $k$ that
$$\lambda(\Omega_N)\geq c\, \varepsilon_N^{n-3} N^{n-1}.$$
Putting, for instance,
$$\varepsilon_N = N^{(1-n)/(n-5/2)},$$
we see that $\Gamma(\Omega_N) \to \infty$, and the desired counterexample is constructed for $n>3$.

\section{Capacitary improvement of the Faber-Krahn inequality}

We state and prove the main result of this section. Here $\Omega$ is an open subset of an arbitrary $n$-dimensional Riemannian manifold.

\smallskip

{\bf Theorem 1.} {\it
Let ${\cal R}>0$, $u\in C_0^\infty(\Omega)$,  and $N_t = \{x\in \Omega: \, |u(x)|\geq t \}$. If $n>2$, then
\begin{eqnarray}\label{k10}
&& \Bigl (\frac{j_{(n-2)/2}}{{\cal R}}\Bigr)^2m_n(B_{\cal R})  \int_0^\infty \Bigl( \frac{{\rm cap}(N_t; \Omega)}{{\rm cap}(B_{\cal R}) + {\rm cap}(N_t; \Omega)} \Bigr)^{\frac{n}{n-2}} d(t^2)\nonumber\\
&& \leq \|\nabla u\|_{L_2(\Omega)}^2,
 \end{eqnarray}
 where $j_\nu$ is the first positive root of the Bessel function $J_\nu$. If $n=2$, then}
 \begin{equation}\label{k11}
 \pi j_0^2 \int_0^\infty {\rm exp}\Bigl(\frac{-4\pi}{{\rm cap}(N_t; \Omega)} \Bigr) d(t^2) \leq \|\nabla u\|_{L_2(\Omega)}^2.
 \end{equation}

\smallskip

{\bf Proof.} Let $w$ be an arbitrary absolutely continuous function on $(0,{\cal R}]$,  such that $w({\cal R}) =0$. The  inequality
 \begin{equation}\label{k12}
 \Bigl (\frac{j_{(n-2)/2}}{{\cal R}}\Bigr)^2\int_0^{\cal R} w(\rho)^2 \rho^{n-1}d\rho \leq \int_0^{\cal R}  w'(\rho)^2 \rho^{n-1}d\rho,
 \end{equation}
where $n>2$, is equivalent to the fact that the first eigenvalue of the Dirichlet-Laplace operator in the unit ball $B$ equals $j_{(n-2)/2}^2$. Similarly, with $n=2$ the inequality
\begin{equation}\label{k13}
 \Bigl (\frac{j_{0}}{{\cal R}}\Bigr)^2\int_0^{\cal R} w(\rho)^2 \rho\, d\rho \leq \int_0^{\cal R} w'(\rho)^2 \rho \, d\rho
 \end{equation}
is associated.

\smallskip

In the case $n>2$, we introduce the new variables
$$\psi = \frac{\rho^{2-n} - {\cal R}^{2-n}}{(n-2)|S^{n-1}|}, \qquad t(\psi) = w(\rho(\psi)),$$
and write (\ref{k12}) in the form
\begin{eqnarray}\label{k14}
&&\bigl(|S^{n-1}| j_{(n-2)/2} {\cal R}^{-1} \bigr)^2 \int_0^\infty \frac{t(\psi)^2 d\psi}{\bigl((n-2)|S^{n-1}|\,
\psi + {\cal R}^{2-n}\bigr)^{2(n-1)/(n-2)}}\nonumber\\
&& \leq \int_0^\infty t'(\psi)^2 d\psi.
 \end{eqnarray}
 Similarly, for $n=2$, putting
$$\psi = (2\pi)^{-1}\log \frac{{\cal R}}{\rho}, \qquad t(\psi) = w(\rho(\psi)),$$
we write (\ref{k13}) as
 \begin{equation}\label{k15}
(2\pi j_0)^2 \int_0^\infty t(\psi)^2 {\rm exp} ( -4\pi\psi) \, d\psi \leq \int_0^\infty t'(\psi)^2 d\psi.
\end{equation}

\noindent
Note that the function $t$ in (\ref{k14}) and (\ref{k15}) is subject to the boundary condition $t(0)=0$. We write (\ref{k14}) and (\ref{k15}) as
 \begin{eqnarray}\label{k16}
&& n^{-1} |S^{n-1}| \Bigl(\frac{j_{(n-2)/2}}{{\cal R}}\Bigr)^2 \int_0^\infty \frac{dt(\psi)^2 }{\bigl((n-2)|S^{n-1}|\, \psi + {\cal R}^{2-n}
\bigr)^{n/(n-2)}}\nonumber\\
&& \leq \int_0^\infty t'(\psi)^2 d\psi
\end{eqnarray}
and
\begin{equation}\label{k17}
\pi j_0^2 \int_0^\infty  {\rm exp} ( -4\pi\psi) \, dt(\psi)^2 \leq \int_0^\infty t'(\psi)^2 d\psi.
\end{equation}

\smallskip

Now, as in Sect. 2.2.1 \cite{M10}, we introduce the function
\begin{equation}\label{k17a}
\psi(t) = \int_0^t \frac{d\tau}{\displaystyle{\int_{|u| =\tau} |\nabla u|\, d H_{n-1}}},
\end{equation}
as well as its inverse $\psi \to t(\psi)$, 
and  replace the integral in the right-hand side of (\ref{k16}) and (\ref{k17}) by $\|\nabla u\|_{L_2(\Omega)}^2$.
It remains to note that

\begin{equation}\label{k17b}
\psi \leq \bigl( {\rm cap} (N_{t(\psi)}; \Omega) \bigr)^{-1}
\end{equation}
by Lemma 2.2.2/1 \cite{M10}.

\smallskip

Let us use the area minimizing function of $\Omega$:
\begin{equation}\label{k17c}
\lambda(v) = \inf H_{n-1}(\partial g),
\end{equation}
where the infimum is extended over all sets $g$ with smooth boudaries  and compact closures   $\overline{g}\subset \Omega$,
subject to the inequality $m_n(g) \geq v.$ 
This and related geometrical characterizations of $\Omega$ proved to be useful in the theory of Sobolev spaces and elliptic equations,  see \cite{M1}, \cite{M2a}, \cite{BM}, \cite{M8}. The function $\lambda$ 
appears in the lower estimate of the capacity
$${\rm cap}(F; \Omega) \geq \Bigl( \int_{m_n(F)}^{m_n(\Omega)} \frac{dv}{\lambda(v)^2} \Bigr)^{-1}$$
(see Corollary 2.2.3/2 \cite{M10}).  Therefore, (\ref{k10}),   (\ref{k11}), and the identity
$${\rm cap}(B_{\cal R}) = (n-2)|S^{n-1}|{\cal R}^{n-2}$$
 lead to the following Lorentz-type estimates.

\medskip

{\bf Corollary 1.} {\it If $n>2$ and ${\cal R} >0$, then, for all $u\in C_0^\infty(\Omega)$,
\begin{eqnarray}\label{k18}
&& \Bigl (\frac{j_{(n-2)/2}}{{\cal R}}\Bigr)^2m_n(B_{\cal R})  \int_0^\infty \Bigl( {\rm cap}(B_{\cal R}) \int_{m_n(N_t)}^{m_n(\Omega)} \frac{dv}{\lambda(v)^2} +1\Bigr)^{\frac{n}{2-n}} d(t^2)\nonumber\\
&&\leq  \|\nabla u\|_{L_2(\Omega)}^2.
\end{eqnarray}
If $n=2$, then, for all}  $u\in C_0^\infty(\Omega)$,
\begin{equation}\label{k19}
\pi j_0^2 \int_0^\infty  {\rm exp}\Bigl ( -4\pi \int_{m_n(N_t)}^{m_n(\Omega)} \frac{dv}{\lambda(v)^2} \Bigr) d(t^2) \leq  \|\nabla u\|_{L_2(\Omega)}^2.
\end{equation}

\medskip

{\bf Remark 1.} Since
\begin{equation}\label{k19a}
\lambda(v) \geq n^{\frac{n-1}{n}} |S^{n-1}|^{\frac{1}{n}} \, v^{\frac{n-1}{n}}
\end{equation}
by the classical isoperimetric inequality for $\Bbb{R}^n$, the estimates (\ref{k18}) and (\ref{k19}) imply the Faber-Krahn property
$$\Lambda(\Omega) \geq \Bigl(\frac{j_{(n-2)/2}}{{\cal R}}\Bigr)^2 $$

\noindent
for any $n$-dimensional Euclidean domain $\Omega$ with $m_n(\Omega) = n^{-1} |S^{n-1}| \, {\cal R}^n$. $\square$

\smallskip

 Theorem 1 is a very special case of the following general assertion.

\medskip

{\bf Theorem 2.}
{\it Let ${\bf M}$ be a decreasing nonnegative function on $[0, \infty)$ and let $q>0$ and $p\geq 1$. Suppose that for all absolutely continuous functions $\psi \to t(\psi)$ on $[0, \infty)$, the inequality
\begin{equation}\label{k20}
\Bigl( -\int_0^\infty |t(\psi)|^q d{\bf M}(\psi) \Bigr)^{1/q} \leq \Bigl( \int_0^\infty |t'(\psi)|^p d\psi  \Bigr)^{1/p}
\end{equation}
holds. Then, for all $u\in C_0^\infty(\Omega)$,
\begin{equation}\label{k21}
\Bigl( \int_0^\infty {\bf M}\Bigl( \bigl({{\rm cap}_p(N_t; \Omega)\bigr)^{1/(1-p)}} \Bigr) d(t^q)\Bigr)^{1/q} \leq \|\nabla u\|_{L_p(\Omega)},
\end{equation}
where ${\rm cap}_p$ is the $p$-capacity defined by}
\begin{equation}\label{k22}
{\rm cap}_p(F;\Omega) = \inf \Bigl\{ \int_\Omega |\nabla u|^p dx: \, u\in C_0^\infty(\Omega), \, u\geq 1 \,\, {\rm on}\,\, F \Bigr\}.
\end{equation}

\smallskip

{\bf Proof.} The role of the function $\psi$ given by (\ref{k17a})  is played in the present proof by
\begin{equation}\label{k23}
\psi(t) = \int_0^t \frac{d\tau}{\displaystyle{\Bigl(\int_{|u| =\tau} |\nabla u|^{p-1}\, d H_{n-1}\Bigr)^{1/(p-1)}}}.
\end{equation}
We write the left-hand side of (\ref{k20}) in the form
$$\Bigl( \int_0^\infty {\bf M}(\psi)\, d(t(\psi))^q \Bigr)^{1/q}$$
and use the monotonicity of ${\bf M}$ and the inequality
\begin{equation}\label{k24}
\psi \leq \bigl({\rm cap}_p(N_{t(\psi)}; \Omega)\bigr)^{1/(1-p)}
\end{equation}
proved in Lemma 2.2.2/1 \cite{M10}. It remains to apply (\ref{k20}) and the identity
\begin{equation}\label{k28a}
\int_0^\infty  |f'(\psi )|^p d\psi =\int _\Omega |\nabla u|^p dx
\end{equation}
found in Lemma 2.3.1 \cite{M10}. $\square$

\smallskip

Using the area minimizing function $\lambda$ defined by (\ref{k17c}) and the estimate
\begin{equation}\label{k25}
{\rm cap}_p(F;\Omega) \geq \Bigl(\int_{m_n(F)} ^{m_n(\Omega)} \frac{dv}{\lambda (v)^{p/(p-1)}}\Bigr)^{1-p}
\end{equation}
(see Corollary 2.2.3/2 \cite{M10}), we obtain from Theorem 2

\smallskip

{\bf Corollary 2.} {\it Let $\mu$, $p$, and $q$ be the same as in Theorem $2$ and let $(\ref{k20})$   hold. Then
\begin{equation}\label{k26}
\Bigl( \int_0^\infty {\bf M} \Bigl(\int_{m_n(N_t)} ^{m_n(\Omega)} \frac{dv}{\lambda (v)^{p/(p-1)}}\Bigr) d(t^q)\Bigr)^{1/q} \leq \|\nabla u\|_{L_p(\Omega)}
\end{equation}
for all} $u\in C_0^\infty(\Omega)$.

\smallskip

Clearly, (\ref{k26}) is a generalization  of the estimates (\ref{k18}) and (\ref{k19}) which were obtained for $p=2$ with  a
particular choice of $\mu$.
Another obvious remark is that (\ref{k20}), where ${\bf M}$ is defined on the interval $0<t<m_n(\Omega)$ by
$$ {\bf M} \Bigl(\int_{t} ^{m_n(\Omega)} \frac{dv}{\lambda (v)^{p/(p-1)}}\Bigr) = \Lambda_{p,q} t$$
with a constant $\Lambda_{p,q}$ depending on $m_n(\Omega)$, implies the inequality
\begin{equation}\label{k31}
\Lambda_{p,q}^{1/q} \|u\|_{L_q(\Omega)} \leq \|\nabla u\|_{L_p(\Omega)}
\end{equation}
for all $u\in C_0^\infty (\Omega)$.

\section{Criterion for an upper estimate of a difference seminorm (the case $p=1$)}

Let us consider the seminorm
\begin{equation}\label{k35t}
\langle u \rangle_{q,\mu} = \Bigl(\int_{\Omega} \int_{\Omega} |u(x) -u(y)|^q \mu(dx, dy)\Bigr)^{1/q},
\end{equation}
where $\Omega$ is an open subset of a Riemannian manifold and $\mu$ is a non-negative measure on
$\Omega\times\Omega$,
locally finite  outside the diagonal $\{(x,y):\;x=y\}$. By definition, the product
$0\cdot \infty $ equals zero.

\smallskip

In this section, first,  we characterize both $\mu$ and $\Omega$ subject to the inequality
\begin{equation}\label{k32}
\langle u \rangle_{q,\mu} \leq C\, \|\nabla u\|_{L_1(\Omega)},
\end{equation}
where $q\geq 1$ and $u$ is an arbitrary function in $C^\infty (\Omega )$. We show that (\ref{k32}) is equivalent
to a somewhat unusual relative isoperimetric inequality.

\medskip

{\bf Theorem 3.} {\it Inequality $(\ref{k32})$  holds for all $u \in C^\infty (\Omega )$ with $q \geq 1$ if and only if
for any open subset $g$ of $\Omega $, such that   $\Omega \cap \partial g\ $ is smooth, the inequality
\begin{equation}\label{k33}
\bigl( \mu(g, \Omega\backslash  \overline{g})+\mu(\Omega\backslash  \overline{g}, {g})\bigr)^{1/q} \leq C H _{n-1}
( \Omega \cap   \partial g )
\end{equation}
holds with the same value of $C$ as in $(\ref{k32})$. In particular, a constant $C$ in $(\ref{k32})$ exists if and only if
}
\begin{equation}\label{k35a}
\sup\limits_{\{g\}} \frac{\mu(g, \Omega\backslash \overline{g} ) ^{1/q}}{H_{n-1} (\Omega\cap \partial g)} <\infty.
\end{equation}

\smallskip

{\bf Proof.} {\it Sufficiency.}  Denote by $u_+$ and $u_-$ the positive and negative parts of $u$, so that $u=u_{+} - u_{-}$.
We notice that
\begin{equation}\label{k34}
\langle u \rangle_{q,\mu} \leq \langle u _{+}\rangle_{q,\mu} + \langle u _{-}\rangle_{q,\mu}
\end{equation}
and
\begin{equation}\label{k35}
\int _\Omega |\nabla u |dx =\int _\Omega |\nabla u _{+}|dx + \int _\Omega |\nabla u _{-}|dx.
\end{equation}
First, we obtain (\ref{k32}) separately for for $u=u_{+}$ and $u=u_{-}$. Let $a>b$ and let $\chi _t(a, b)=1$ if $a>t>  b$
and $\chi _t(a, b)=0$ otherwise.

\smallskip

Clearly,
\begin{eqnarray*}
\langle u \rangle_{q,\mu}&=&\left ( \int _\Omega  \int _\Omega  \left | \int _{u(x)} ^{u(y)}dt \right |^q   \mu (dx,dy)\right  )^{1/q}\\
&=&\left (\int _\Omega\int _\Omega  \left | \int _{0} ^{\infty } \big ( \chi _t(u(x), u(y)) +\chi _t (u(y), u(x)) \big ) dt \right |^q  \mu (dx,dy)
\right ) ^{1/q}.
\end{eqnarray*}
By Minkowski's inequalitiy,
\begin{eqnarray*}
\langle u \rangle_{q,\mu}&\leq &
\int _{0} ^{\infty }\left ( \int _\Omega\int _\Omega  \big ( \chi _t(u(x), u(y)) +\chi _t (u(y), u(x)) \big ) ^q  \mu (dx,dy)\right ) ^{1/q}dt\\
&=&\int _{0} ^{\infty }\left ( \int _\Omega\int _\Omega  \big ( \chi _t(u(x), u(y)) +\chi _t (u(y), u(x)) \big )\mu (dx,dy)  \right ) ^{1/q} dt\\
&=&\int _{0} ^{\infty }\big ( \mu (M_t, \Omega \backslash N_t)+ \mu (\Omega \backslash N_t, M_t )\big )^{1/q}dt,
\end{eqnarray*}
where $M_t=\{x \in \Omega  :u(x) >t \}$ and $N_t=\{x \in \Omega  :u(x) \geq t \}$.

\smallskip

By (\ref{k33}) and the co-area formula, the last integral does not exceed
$$
C\int _{0} ^{\infty } H_{n-1} \big (\{x \in \Omega  :u(x) >t \} \big )dt=C\int _\Omega |\nabla u(x)|dx.
$$
Therefore,
$$
\langle u _{\pm }\rangle_{q, \mu}\leq C\int _\Omega |\nabla u_{\pm }(x)|dx
$$
and the reference to  (\ref{k34}) and  (\ref{k35}) completes the proof of sufficiency.

\smallskip

{\it Necessity.} Let $\{ w_m \}$ be the sequence of locally Lipschitz functions in $\Omega$ constructed in 
 Lemma  3.2.2  \cite{M10} with the following properties:
\begin{enumerate}
 \item $w_m(x) =0$ in $\Omega\backslash g$,
 \item $w_m(x)\in [0, 1]$ in $\Omega$,
 \item for any compactum $K\subset g$ there exists an integer $N(e)$ such that $w_m(x) =1$ for $x\in K$ and $m\geq N(e)$,
\item the limit relation holds
$$
 \mathop{\hbox {lim sup}}_{m\to \infty} \int_\Omega |\nabla w_m(x)|dx = H_{n-1}(\Omega\cap \partial g).
$$
\end{enumerate}
 By  Theorem 1.1.5/1 \cite{M10}, the inequality  (\ref{k32}) holds for all locally Lipschitz functions. Therefore,
\begin{equation}\label{k36t} 
\langle w_m \rangle_{q,\mu}\leq C\, \|\nabla w_m \|_{L_1(\Omega )} 
 \end{equation}
and due to property 4, 
\begin{equation}\label{k36} 
\mathop{\hbox {lim sup}}_{m\to \infty} \langle w_m \rangle_{q,\mu}\leq CH_{n-1}(\Omega \cap \partial g).
\end{equation}
On the other hand,
\begin{eqnarray*}
\langle w_m \rangle_{q,\mu}^q&= &
\int _{x\in g}\int _{y\in \Omega \backslash g} w_m(x) ^q  \mu (dx,dy)\\
&+& \int _{x\in \Omega \backslash g}\int _{y \in g} w_m(y) ^q  \mu (dx,dy)+ \int _{g}\int _{g}|w_m(x)-w_m(y)|^q \mu (dx,dy)
\end{eqnarray*}
which implies
$$
\langle w_m \rangle_{q,\mu}^q \geq 
\int _{ g} w_m(x) ^q  \mu (dx, \Omega \backslash \overline{g})+\int _{ g} w_m(y) ^q  \mu (\Omega \backslash \overline{g}, dy ).
$$
This, along with  property 3,  leads to
$$
\mathop{\hbox {lim inf}}_{m\to \infty}\, \langle w_m \rangle_{q,\mu}^q   \geq  \mu (g, \Omega \backslash \overline{g})+
\mu (\Omega \backslash g, \overline{g} ).
$$
Combining this relation with (\ref{k36t}) and (\ref{k36}), we arrive at (\ref{k33}). $\square$

\medskip

{\bf Corollary 3} (One-dimensional case). {\it Let
$$
\Omega = (\alpha, \beta), \quad where \quad -\infty\leq \alpha <\beta \leq \infty.
$$
The inequality
\begin{equation}\label{k37}
\Bigl(\int_{\Omega} \int_{\Omega} |u(x) -u(y)|^q \mu(dx, dy)\Bigr)^{1/q} \leq C \int_\Omega |u'(x)| dx
\end{equation}

\noindent
with $q\geq 1$ holds for all $u\in C^\infty(\Omega)$ if and only if
\begin{equation}\label{k38}
\bigl( \mu(I, \Omega\backslash \overline{I})+\mu(\Omega\backslash  \overline{I}, I)\bigr)^{1/q} \leq 2C
\end{equation}
for all open intervals $I$ such that  $\overline I \subset \Omega $, and
\begin{equation}\label{k39}
\bigl( \mu(I, \Omega\backslash \overline{I})+\mu(\Omega\backslash  \overline{I}, I)\bigr)^{1/q} \leq C
\end{equation}
for all intervals $I \subset \Omega $, such that  $\overline I $ contains one of the end points of $\Omega $.

\smallskip

In particular, a constant in $(\ref{k37})$ exists if and only if
}
$$\sup\limits_{\{I\}} \mu (I, \Omega\backslash \overline{I}) <\infty.$$

\smallskip
{\bf Proof.} {\it Necessity} follows directly from (\ref{k33}) by setting $g=I$. Let us check the sufficiency of (\ref{k38}).
Represent an arbitrary open set $g$ of $\Omega $ as the union of non-overlapping open intervals $I_k$. Then by
(\ref{k38}) and (\ref{k39})
\begin{eqnarray*}
& &\bigl( \mu(g, \Omega\backslash \overline{g})+\mu(\Omega\backslash  \overline{g}, g)\bigr)^{1/q}
=\left (\sum_{k}\bigl( \mu(I_k, \Omega\backslash \overline{g})+\mu(\Omega\backslash  \overline{g}, I_k)\bigr) \right )^{1/q} \\
& &\leq\sum_{k}\bigl( \mu(I_k, \Omega\backslash \overline{g})+\mu(\Omega\backslash  \overline{g}, I_k)\bigr)^{1/q}\leq C\sum_{k}
H_0(\Omega \cap \partial I_k)
\end{eqnarray*}
which is the same as (\ref{k33}). The result follows from Theorem 3.$\square$

\medskip

{\bf Remark 2.} Suppose that the class of admissible functions in Theorem 3 is diminished by the requirement 
$u=0$ {\it in a neighbourhood of a closed subset $F$ of} $\overline \Omega $. Then the same proof leads to the same
criterion (\ref{k33}) with the only difference that the admissible sets $g$ should be at a positive distance from $F$.
For the example $F=\partial \Omega $, i.e. for the inequality (\ref{k32}) with any $u \in C^\infty _0(\Omega )$,
the necessary and sufficient condition (\ref{k33}) becomes the  isoperimetric inequality
\begin{equation}\label{k40}
\bigl( \mu(g, \Omega\backslash \overline{g})+\mu(\Omega\backslash \overline{g}, g)\bigr)^{1/q} \leq C H_{n-1}(\partial g)
\end{equation}
for all open sets $g$ with smooth boundary and compact closure  $\overline{g} \subset \Omega $. If, in particular, in Corollary 3,
the criterion of the validity of (\ref{k37}) for all $u \in C^\infty _0(\Omega )$ is the inequality (\ref{k38})   for
every interval $I$, $\overline I \subset \Omega $. In the case $u=0$ near one of the end points $\Omega =(\alpha ,\beta )$,
one should require both (\ref{k38}) and (\ref{k39}) but the intervals $I$ should be at a positive distance from that end point.

\smallskip

Needless to say, the condition (\ref{k33}) is simplified as follows for a symmetric measure $\mu $, i.e. under the assumption
$\mu({\cal E}, {\cal F}) = \mu( {\cal F}, {\cal E})$:
$$
\mu(g, \Omega\backslash \overline{g})^{1/q}\leq 2^{-1/q}CH_{n-1}(\Omega \cap \partial g)
$$
for the same open sets $g$ as in Theorem 3.

\medskip

{\bf Remark 3.} The integration domain $\Omega\times \Omega$ 
in (\ref{k35t}) excludes inequalities for integrals taken over $\partial\Omega$. This can be easily avoided assuming additionally that $\mu$ is defined on compact subsets of $\overline{\Omega}\times \overline{\Omega}$ and that $u\in C(\overline{\Omega})\cap C^\infty(\Omega)$. Then, with the same proof, one obtains the corresponding criterion, similar to (\ref{k33}):
$$\bigl( \mu(\overline{g}, \, \overline{\Omega}\backslash  \overline{g})+\mu(\overline{\Omega}\backslash  \overline{g}, \, \overline{g})\bigr)^{1/q} \leq C\,  H _{n-1}
( \Omega \cap   \partial g ).
$$
As an application, consider the inequality
\begin{equation}\label{k88}
\int_{\partial\Omega}\int_{\partial\Omega} |u(x) -u(y)|\, H_{n-1}(dx)\, H_{n-1}(dy) \leq C \int_\Omega |\nabla u|\, dx
\end{equation}
which holds if and only if
\begin{equation}\label{k89}
H _{n-1}(\partial \Omega \cap   \partial g )\, H _{n-1}( \partial\Omega \backslash   \partial g ) \leq 2^{-1} C\, H _{n-1}( \Omega \cap   \partial g )
\end{equation}
for the same sets $g$ as in Theorem 3.

\smallskip

By Corollary 6.4.4/3 \cite{M10}, which appeared first in \cite{BM}, 

(i) If $\Omega$ is the unit ball in $\Bbb{R}^3$, then
$$4\pi H _{2}( \Omega \cap   \partial g ) \geq H _{2}( \partial\Omega \cap   \partial g )\, H _{2}( \partial\Omega \backslash   \partial g )$$
and

(ii) If $\Omega$ is the unit disk on the plane, then
$$H _{1}( \Omega \cap   \partial g ) \geq 2\sin \bigl(\frac{1}{2} H _{1}( \partial\Omega \cap   \partial g ) \bigr).
$$
Moreover, the last two inequalities are sharp. Hence, the inequality (\ref{k88}) holds with the best constant $C = 8\pi$ if $\Omega = B$. In the case (ii),
$$
H _{1}( \Omega \cap   \partial g ) \geq 2^{-1} \min\limits_{0\leq\varphi \leq \pi} \frac{\sin\varphi}{\varphi (\pi-\varphi)} H _{1}( \partial\Omega \cap   \partial g )H _{1}( \partial\Omega \backslash   \partial g ).$$
Since the last minimum equals $\pi^{-1}$, it follows that the best value of $C$ in the inequality (\ref{k88}) for the unit disk is $4\pi$. $\square$

\medskip

We can simplify the criterion (\ref{k33}) for $\Omega = \Bbb{R}^n$, replacing arbitrary sets $g$ by arbitrary balls $B(x,\rho)$ similarly to Theorem 1.4.2/2 \cite{M10}, where the norm 
$$\|u\|_{L_q(\mu)} = \Bigl( \int_{ \Bbb{R}^n} |u|^q d\mu \Bigr)^{1/q}$$
is treated in place of $\langle u \rangle_{q,\mu}$. Unfortunately, the best constant in the sufficiency part will be lost.

\smallskip

{\bf Corollary 4.} (i) {\it If $q\geq 1$ and 
\begin{equation}\label{k90}
\sup\limits_{x\in  \Bbb{R}^n, \rho >0} \rho^{(1-n)q}\bigl(\, \mu(B(x,\rho), \Bbb{R}^n \backslash B(x,\rho)) + 
\mu( \Bbb{R}^n \backslash B(x,\rho), B(x,\rho))\, \bigr) <\infty,
\end{equation}
then the inequality
\begin{equation}\label{k91}
\Bigl(  \int_{ \Bbb{R}^n} \int_{ \Bbb{R}^n} |u(x) -u(y)|^q \mu(dx,dy)\Bigr)^{1/q} \leq C\, \|\nabla u\|_{L_1( \Bbb{R}^n)}
\end{equation}
holds for all $u\in C^\infty( \Bbb{R}^n)$    and
\begin{equation}\label{k92}
C^q \leq c^q \sup\limits_{x\in  \Bbb{R}^n, \rho >0} \rho^{(1-n)q}\bigl(\, \mu(B(x,\rho), \Bbb{R}^n \backslash B(x,\rho) )+ 
\mu( \Bbb{R}^n \backslash B(x,\rho), B(x,\rho)) \, \bigr),
\end{equation}
where $c$ depends only on $n$.}

\smallskip

(ii) {\it If $(\ref{k91})$ holds for all $u\in C^\infty( \Bbb{R}^n)$, then}
$$
C^q \geq |S^{n-1}|^{-q} \sup\limits_{x\in  \Bbb{R}^n, \rho >0} \rho^{(1-n)q}\bigl(\, \mu(B(x,\rho), \Bbb{R}^n \backslash B(x,\rho)) + 
\mu( \Bbb{R}^n \backslash B(x,\rho), B(x,\rho)) \, \bigr).
$$

\smallskip

{\bf Proof.} Let $g$ be an arbitrary open set in $\Bbb{R}^n$ with  smooth boundary and let $\{B(x_j,\rho_j)\}$ be the Gustin covering of $g$ subject to
\begin{equation}\label{k93}
\sum_j \rho_j^{n-1} \leq c\, H_{n-1} (\partial g),
\end{equation}
where $c$ depends only on $n$ (see Theorem 1.2.2/2 \cite{M10}).  Then
\begin{eqnarray*}
 \mu(g, \Bbb{R}^n\backslash g) &\leq & \sum_j \mu( B(x_j,\rho_j), \Bbb{R}^n \backslash g)\\
 &\leq &\Bigl( \sum_j \mu( B(x_j,\rho_j), \Bbb{R}^n \backslash g)^{1/q}\Bigr)^q\\
 &\leq &\Bigl( \sum_j \mu( B(x_j,\rho_j), \Bbb{R}^n \backslash B(x_j,\rho_j))^{1/q} \Bigr)^q\\
&\leq & (c\, B)^q \bigl(\sum_j \rho_j^{n-1}\bigr)^q,
\end{eqnarray*}
where $B$ is the value of the supremum in (\ref{k90}). This and (\ref{k93}) imply
$$\mu(g, B(x_j,\rho_j) \leq (c\, B\, H_{n-1}(\partial g))^q.$$
Similarly,
$$\mu(\Bbb{R}^n \backslash g, g) \leq (c\, B\, H_{n-1}(\partial g))^q$$

\noindent
and the result follows from Theorem 3.

The assertion (ii) stems from (\ref{k33}) by setting $g= B(x,\rho)$. $\square$

\section{Criterion for an upper estimate of a difference norm (the case $p>1$)}

Now we deal with the inequality
\begin{equation}\label{k41}
\langle u \rangle_{q,\mu} \leq C\, \|\nabla u\|_{L_p(\Omega)},
\end{equation}
where $q >p> 1$, and show that it is equivalent to a certain isocapacitary inequality.

\smallskip

The capacity to appear in the present context is defined as follows. Let $F_1$ and $F_2$ be non-overlapping
subsets of $\Omega $, closed in $\Omega $. The $p$-capacity of the pair $(F_1, F_2)$ with respect to $\Omega $
is given  by
$$
 {\rm cap} _p(F_1, F_2; \Omega)=\inf _{\{u\}}\int _\Omega |\nabla u(x)|^p dx,
$$
where $\{ u \}$ is the set of all $u \in C^\infty (\Omega )$, such that $u \geq 1$ on $F_1$ and $u\leq 0$ on $F_2$.

\smallskip

Obviously, this capacity does not change if $F_1$ and $F_2$ change places. Furthermore, if $F$ is a closed set in ${\Bbb R}^n$
and $F\subset G$, where $G$ is an open set, such that  $\overline G \subset \Omega $,  then
$ {\rm cap} _p(F, \Omega \backslash G; \Omega )$ coincides with the $p$-capacity ${\rm cap} _p(F; G )$ defined in
(\ref{k22}).

\medskip

{\bf Theorem 4}. {\it Inequality $(\ref{k41})$ with $p\in (1,q)$  holds for all $u \in C^\infty (\Omega )$ if and only if
for any pair $(F_1, F_2)$ of non-overlapping sets, closed in $\Omega $,
\begin{equation}\label{k42}
\mu (F_1, F_2)^{p/q}\leq  B\;{\rm cap} _p(F_1, F_2; \Omega),
\end{equation}
where $B$ depends only on $p$ and $q$. In the sufficiency part we may assume that $F_1$ and $F_2$ are sets with smooth
$\Omega \cap \partial F_i$.}

\smallskip
In the proof of this theorem, we use the  inequality
\begin{equation}\label{k42b}
\Bigl( \int_ {{\Bbb R}_+}  |f(\psi)|^q\psi^{-1-q/p'} d\psi \Bigr)^{1/q}  \leq c\,  \|f'\|_{L_p({\Bbb R}_+)}
\end{equation}
due to Bliss  \cite{Bl} and the inequality
\begin{equation}\label{k42a}
\Bigl ( \int_ {{\Bbb R}_+} \int_ {{\Bbb R}_+} {|f(\psi )-f(\phi )|^q \over |\psi -\phi |^
{2+q/p'}}d\phi d\psi \Bigr )^{1/q}\leq c \, \|f' \|_{L_p({\Bbb R}_+)},
\end{equation}
where $q>p>1, p'=p/(p-1)$ and $f$ is an arbitrary absolutely continuous function on ${\overline {\Bbb R}_+}$.

\smallskip

A short argument leading to (\ref{k42a}) is as follows. 
Clearly, (\ref{k42a}) results   from the same inequality  with  ${\Bbb R}$  in place of ${\Bbb R}_+$,
which follows, in its turn, from the estimate
\begin{equation}\label{k43a}
\|f \|_{B^{1-(q-p)/pq}_{q}({\Bbb R})}\leq c\,  \|f \|_{W^1_p({\Bbb R})}
\end{equation}
by dilation with a coefficient $\lambda $ and the limit passage as $\lambda  \rightarrow 0_+$. (The standard notations $B$ and $W$ for Besov and Sobolev spaces with non-homogeneous norms is used in (\ref{k42a}).)
In order to obtain
(\ref{k43a}), we recall the well-known Sobolev type inequality
$$
\|h \|_{L_{p'}({\Bbb R}_+)}\leq c \, \|h \|_{B^{(q-p)/pq}_{q'}({\Bbb R})}
$$
(see Theorem $4'$, Sect. 5.1  \cite{St}) and put $h=(-\Delta +1)^{-1/2}f$, which shows  that
\begin{equation}\label{k44a}
\|f \|_{W^{-1}_{p'}({\Bbb R})} \leq c \, \|f \|_{B^{-1+(q-p)/pq}_{q'}({\Bbb R})}.
\end{equation}
By duality, (\ref{k44a}) is equivalent to (\ref{k43a}).

\smallskip

With (\ref{k42a}) at hand, we   return to Theorem 4.

\smallskip

{\bf Proof.} {\it Sufficiency.} Arguing as at the beginning of the the proof of Theorem 2, we see that it sufficies to prove (\ref{k41})
  for a non-negative $u$. By the definition of the Lebesque integral
$$
\int _\Omega u d\nu =\int _{{\Bbb R}_+} \nu (N_\tau )d\tau =\int _{{\Bbb R}_+} \nu (M_\tau )d\tau,
$$
where $\nu $ is a measure, and therefore
\begin{equation}\label{k45a}
\int _\Omega P(u)d\nu =\int _{{\Bbb R}_+}\nu (N_\tau )dP(\tau ),
\end{equation}
where $P$ is a non-decreasing function on ${\Bbb R}_+$. Putting  here $u=1/v$ and  $Q(\tau) =P(\tau ^{-1})$,
we deduce
\begin{equation}\label{k46a}
\int _\Omega Q(u)d\nu =- \int _{{\Bbb R}_+}\nu (\Omega \backslash M_\tau )dQ(\tau ),
\end{equation}
where $Q$ is non-increasing. We obtain
\begin{eqnarray*}
\int _\Omega \int _\Omega |u(x)-u(y)|^q \mu(dx, dy)&=&\int _\Omega \int _\Omega (u(x)-u(y))^q_+ \mu(dx, dy)\\
&+&\int _\Omega \int _\Omega (u(y)-u(x))^q_+ \mu(dx, dy)\\
&=&\int _\Omega \int _\Omega (u(x)-u(y))^q_+ (\mu(dx, dy)+\mu(dy, dx)).
\end{eqnarray*}
By (\ref{k45a}) and (\ref{k46a}), the last double integral is equal to
\begin{eqnarray*}
& &q\int _{{\Bbb R}_+} \int _\Omega (t -u(y))_+^{q-1} (\mu(N_\tau , dy)+\mu(dy, N_\tau ))d\tau \\
& &=q(q-1)\int _{{\Bbb R}_+} \int _{{\Bbb R}_+} (\tau -\sigma )^{q-2}_+\bigl ( \mu(N_\tau ,\Omega \backslash M_\sigma )
+\mu(\Omega \backslash M_\sigma, N_\tau )\bigr )d\tau d\sigma .
\end{eqnarray*}
Now, (\ref{k42}) implies
$$
\langle u \rangle_{q,\mu}^q \leq 2q(q-1)B\int _{{\Bbb R}_+}\int _{{\Bbb R}_+}(\tau -\sigma )^{q-2}_+ {\rm cap}_p
(N_\tau ,\Omega \backslash M_\sigma ; \Omega )d\tau d\sigma
$$
and using the function $\psi \rightarrow t(\psi )$, inverse of (\ref{k23}), we arrive at the inequality
\begin{eqnarray*}
& &\hspace{-15mm}\langle u \rangle_{q,\mu}^q \leq 2q(q-1)B^{q/p}\\
& &\times \int _{{\Bbb R}_+}\int _0 ^\psi (t(\psi ) -t(\phi ) )^{q-2}  \bigl ({\rm cap}
(N_{t(\psi )} ,\Omega \backslash M_{t(\phi )} ; \Omega )^{q/p}t'(\phi )t'(\psi )d\phi  d\psi .
\end{eqnarray*}
By Lemma 2.2.2/1 \cite{M10}, for $\psi >\phi $
$$
{\rm cap}
(N_{t(\psi )} ,\Omega \backslash M_{t(\phi )} ; \Omega ) \leq (\psi -\phi )^{1-p}
$$
and  therefore,
\begin{equation}\label{k47a}
\langle u \rangle_{q,\mu}^q \leq 2q(q-1)B^{q/p}
\int _{{\Bbb R}_+}\int _0 ^\psi (\psi -\phi )^{-q/p'}(t(\psi ) -t(\phi ) )^{q-2}t'(\phi )t'(\psi )d\phi  d\psi .
\end{equation}
Integrating by parts twice on the right-hand side of (\ref{k47a}), we obtain
$$\langle u \rangle_{q,\mu}^q \leq  2B^{q / p}{q \over p'}\Bigl(\Bigl( {q \over p'} +1\Bigr )
\int_ {{\Bbb R}_+} \!\int_ {0}^\psi  {(t(\psi )-t(\phi ))^q \over (\psi -\phi )^
{2+q/p'}}d\phi d\psi \!+\! \int_ {{\Bbb R}_+} \psi^{-q/p'} t(\psi)^{q-1} t'(\psi) d\psi\Bigr)
$$
$$=B^{q / p}{q \over p'}\Bigl(\Bigl ( {q \over p'} +1\Bigr)
\int_ {{\Bbb R}_+} \int_ {{\Bbb R}_+} {|t(\psi )-t(\phi )|^q \over |\psi -\phi |^{2+q/p'}}d\phi d\psi  + \frac{1}{p'}\int_ {{\Bbb R}_+} t(\psi)^q \psi^{-1-q/p'} d\psi\Bigr).$$
Hence, we deduce from (\ref{k42b}) and (\ref{k42a}) that
\begin{equation}\label{k47b}
\langle u \rangle_{q,\mu}\leq  c\, B^{1/p}\|t' \|_{L_p( {\Bbb R}_+)},
\end{equation}
where $c$ depends only on $p$ and $q$. It remains to refer to (\ref{k28a}).

\smallskip

{\it Necessity}. Let $F_1$ and $F_2$ be subsets of $\Omega $, closed in $\Omega $. We take an arbitrary function
$u \in C^\infty (\Omega )$, such that $u\geq 1$ on $F_1$ and $u \leq 0$ on $F_2$, and put it into (\ref{k41})
$$
\mu (F_1, F_2;\Omega )^{p/q}\leq \left ( \int _{F_1} \int _{F_2} |u(x)-u(y)|^q \mu(dx, dy) \right )^{1/q} \leq
C\int _\Omega |\nabla u|^{p}dx.
$$
It remains to minimize the right-hand side, in order to obtain
$$
\mu (F_1, F_2;\Omega )^{p/q}\leq C\, {\rm cap}_p (F_1, F_2;\Omega ).
$$
The result follows.$\square$

\bigskip

A direct consequence of Theorem 4 and the isocapacitary inequality for ${\rm cap}_p(F; G)$ (see (5) and (6) in Sect. 2.2.3 \cite{M10}) is
the following sufficient condition for (\ref{k41}) formulated in terms of the $n$-dimensional Lebesgue measure:
\begin{equation}\label{k43}
\mu(F, \Omega\backslash  G) \leq c \Bigl( \log \frac{m_n(G)}{m_n(F)} \Bigr)^{q(1-n)/n}, \quad {\rm if} \,\, p=n
\end{equation}
and

\begin{equation}\label{k44}
\mu(F, \Omega\backslash  G) \leq c \bigl | m_n(G)^{(p-n)/n(p-1)} - m_n(F)^{(p-n)/n(p-1)} \bigr |^{1-p}, \quad {\rm if} \,\, p\neq n.
\end{equation}

\smallskip

Choosing two concentric balls situated in $\Omega $ as the sets $F_1$   and $\Omega \backslash F_2$ in (\ref{k42}) and using the explicit
fofmulae for the $p$-capacity of spherical condensers (see (1) and (2) in Sect. 2.2.4 \cite{M10}) we see that the inequalities (\ref{k43})
and (\ref{k44}),  with concentric  balls $F$ and $G$ placed in $\Omega$,  is a necessary condition for (\ref{k41}).

In the one-dimensional case Theorem 4 can be written in a much simpler form.

\smallskip

{\bf Corollary 5.} {\it Let
$$\Omega = (\alpha, \beta), \qquad -\infty\leq \alpha <\beta \leq \infty.$$
The inequality
\begin{equation}\label{k96}
 \Bigl(\int_{\Omega} \int_{\Omega} |u(x) -u(y)|^q \mu(dx, dy)\Bigr)^{1/q} \leq C \Bigl(\int_\Omega |u'(x)|^p dx\Bigr)^{1/p}
 \end{equation}
holds for every $u\in C^\infty(\Omega)$ if and only if, for all pair of intervals $I$ and $J$ of the three types:
$$I = [x-d, x+d] \quad {\rm and} \,\,\, J =  (x-d-r, \, x+d+r),$$
\begin{equation}\label{k96a}
I = (\alpha, x] \quad {\rm and} \,\,\, J = (\alpha, x+r),
 \end{equation}
\begin{equation}\label{k96b}
I = [x, \beta) \quad {\rm and} \,\,\, J = (x-r, \beta),
 \end{equation}
where $d$ and $r$ are positive and $J\subset \Omega$, we have
\begin{equation}\label{k45}
r^{(p-1)/p} \bigl(\mu(I, \Omega\backslash J)\bigr)^{1/q} \leq B,
\end{equation}
where $B$ does not depend on $I$ and $J$.}

\smallskip

{\bf Proof.} The necessity of (\ref{k45}) follows directly from that in Theorem 4 and the inequality
$${\rm cap}_p(I, \Omega\backslash J; \Omega) \leq 2\, r^{1-p}$$
(see Lemma 2.2.2/2 \cite{M10}). 

\smallskip

Let us prove the sufficiency. 
By $G_1$ we mean an open subset of $\Omega $ such that $F_1 \subset G_1$ and
$\overline G_1 \subset \Omega \backslash F_2$. Connected components of $\Omega \backslash F_2$ will be denoted
by $J_k$. Let   $J_k$ contain the closed convex hull $\tilde{J_k}$ of those connected components of $G_1$
which are situated in $J_k$.

\smallskip

Then
$$
\mu (F_1, F_2)^{p/q}\leq  \mu (G_1, F_2)^{p/q}\leq \Bigl (\sum_{k}\mu \bigl ( \tilde{J_k}, \Omega \backslash J_k \bigr )
\Bigr ) ^{p/q} \leq \sum_{k}\mu \bigl ( \tilde{J_k}, \Omega \backslash J_k \bigr ) ^{p/q}
$$
and since by (\ref{k45})
$$
\mu \left ( \tilde{J_k}, \Omega \backslash J_k \right )^{p/q} \leq B^p \bigl ( {\rm dist} \{ I_k, {\Bbb R} \backslash J_k \} \bigr )^{1-p}
$$
we obtain
\begin{equation}\label{k46}
\mu (F_1, F_2)^{p/q}\leq B^p \sum_{k} \bigl ( {\rm dist} \{ I_k, {\Bbb R} \backslash J_k \} \bigr )^{1-p}.
\end{equation}

Consider an arbitrary function $u \in C^\infty (\Omega )$, such that $u=1$ on $G_1$ and $u=0$ on $F_2$. Clearly,
$u=0$ on $\partial J_k$. We have
\begin{equation}\label{k47}
\int _\Omega |u'|^p dx \geq \sum_{k} \int _{J_k} |u'|^p dx \geq \sum_{k} \int _{J_k} |\tilde{u}'_k|^pdx ,
\end{equation}
where $\tilde{u}=u$ on ${J_k}\backslash \tilde{I}_k$, $\tilde{u}_k=1$ on $\tilde{I}_k$,  and $\tilde{u}_k=0$ on $\partial J_k$.
Hence
$$
\int _\Omega |u'|^p dx \geq \sum_{k} \bigl ( {\rm dist} \{ \tilde{I}_k, {\Bbb R} \backslash J_k \} \bigr )^{1-p}.
$$
Comparing this estimate with (\ref{k46}), we arrive at
$$
\int _\Omega |u'|^p dx \geq \mu (F_1, F_2)^{p/q}
$$
and minimizing the integral in the left-hand side over all functions $u$, we obtain (\ref{k45}).$\square$

\medskip

{\bf Remark 4.} It is straightforward but somewhat cumbersome to obtain a more general criterion by replacing the seminorm on the right-hand side of (\ref{k96}) with
\begin{equation}\label{k97}
\Bigl(\int_\Omega |u'(x)|^p \sigma(dx)\Bigr)^{1/p},
\end{equation}
where $\sigma$ is a measure in $\Omega$. In fact, one can replace $\sigma$ by its absolutely continuous part $(d\sigma^*/dx)dx$ and further, roughly speaking, the criterion will follow from Corollary 5 by the change of variable $x\to \xi$, where
$$d\xi = (d\sigma^*/dx)^{1/(1-p)} dx.$$
Restricting myself to this hint, I leave details to the interested reader.

\section{Capacitary sufficient condition in the case $q=p$}

 In the marginal case $q=p$ the condition (\ref{k42}) in Theorem 4, being necessary, is not generally
sufficient. In fact, let $n=1, \Omega={\Bbb R}$, and
$$
\mu (dx, dy)={dxdy \over |x-y|^{p+1}} .
$$
Then as shown in the proof of Corollary 4,  (\ref{k42}) is equivalent to (\ref{k45}),  and (\ref{k45}) holds, since
\begin{eqnarray*}
\mu (I, {\Bbb R}\backslash J)&=&\int _{|t-x|<d}dt\int _{|\tau -x|>d+r}{d\tau \over |t-\tau |^{p+1}}\\
& &\\
&=&\int _{|t|<d}dt\int _{|\tau |>d+r}{d\tau \over |t-\tau |^{p+1}}\leq c\, r^{1-p}
\end{eqnarray*}
and the same estimate holds for $I$ and $J$ defined by (\ref{k96a}) and (\ref{k96b}).

\smallskip

On the other hand, (\ref{k40}) fails, because
$$
\int _{\Bbb R} \int _{\Bbb R} {|u(x)-u(y)|^p \over |x-y|^{p+1}}dxdy=\infty
$$
for every non-constant function $u$. 

\smallskip

In the next theorem we give a sufficient condition for (\ref{k41}) with $q=p>1$ formulated in terms of an isocapacitary
inequality.

\medskip

{\bf Theorem 5.} {\it Given $p\in (1,\infty)$ and a positive, vanishing at infinity, non-increasing absolutely continuous  function $\nu$ on ${\Bbb R}_+$, such that
$$S: =  \sup\limits_{\tau >0} \Bigl(\int_0^\tau |\nu'(\sigma)|^{1/(1-p)}\frac{d\sigma}{\sigma} \Bigr)^{p-1}  \int_\tau ^\infty |\nu'(\sigma)|\frac{d\sigma}{\sigma}<\infty.$$
Suppose that
\begin{equation}\label{k57}
\mu(F_1, F_2) \leq \nu\bigl( ({\rm cap}_p(F_1, F_2; \Omega) )^{1-p}\bigr)
\end{equation}
for all non-overlapping sets $F_1$ and $F_2$ closed in $\Omega$. Assume also that
\begin{equation}\label{k57a}
{\cal K} : = \int_0^\infty |\nu'(\sigma)|\, \sigma^{p-1} d\sigma <\infty.
\end{equation}
Then 
\begin{equation}\label{k58}
\|u\|_{p,\mu} \leq 2^{1/p} p \, \Bigl( \frac{S}{(p-1)^{p-1}} \Bigr) ^{1/p p'} {\cal K}^{1/p} \|\nabla u\| _{L_p(\Omega)}
\end{equation}
for all} $u\in C^\infty(\Omega)$. 

\smallskip

{\bf Proof.} We assume that $\nabla u\in L_p(\Omega)$ and the integral in (\ref{k58}) involving derivatives of $\nu$ is convergent.  Arguing as in the proof of Theorem 4 and using (\ref{k57}) instead of (\ref{k42}), we obtain
\begin{equation}\label{k59}
\langle u\rangle_{p,\mu}^p \leq 2p(p-1) \int_0^\infty\int_\phi^\infty \nu(\psi -\phi) \bigl( t(\psi) - t(\phi)\bigr)^{p-2} t'(\psi) d\psi \, t'(\phi) d\phi.
\end{equation}
Owing to  (\ref{k57a}), we can integrate by parts in the inner integral in (\ref{k59}) and obtain
\begin{eqnarray*}
\langle u\rangle_{p,\mu}^p& \leq &  2p \int_0^\infty \int _\phi^\infty |\nu'(\psi -\phi)| \bigl(t(\psi) -t(\phi)\bigr)^{p-1} d\psi\, t'(\phi) d\phi\\
&=& 2p \int_0^\infty \int _0^\psi |\nu'(\psi -\phi)| \bigl(t(\psi) -t(\phi)\bigr)^{p-1} \, t'(\phi) d\phi \, d\psi.
\end{eqnarray*}
By H\"older's inequality
\begin{equation}\label{k60}
\langle u\rangle_{p,\mu}^p \leq   2p \int_0^\infty {\cal A}(\phi)^{1/p'} {\cal B}(\phi)^{1/p} d\phi,
\end{equation}
where
$${\cal A} = \int_0^\psi \frac{ |\nu'(\psi -\phi)| }{\psi -\phi}\bigl(t(\psi) -t(\phi)\bigr)^{p} d\phi$$
and 
$${\cal B}  = \int_0^\psi { |\nu'(\psi -\phi)| }\, (\psi -\phi)^{p-1}\, |t'(\psi)|^{p} d\phi.$$
Using Theorem 1.3.1/1 \cite{M10} concerning 
a two-weight Hardy inequality, we obtain
$${\cal A} \leq \frac{p^p}{(p-1)^{p-1}}\, S\, {\cal B} $$
which together with (\ref{k60}) gives
$$\langle u\rangle_{p,\mu}^p \leq   2p^p(p-1)^{(1-p)/p'} S^{1/p'}\int_0^\infty
\int_0^\psi { |\nu'(\psi -\phi)| }\, (\psi -\phi)^{p-1}\, |t'(\psi)|^{p} d\phi\, d\psi.$$
Changing the order of integration, we arrive at
$$\langle u\rangle_{p,\mu} \leq   2^{1/p} p\, \bigl( (p-1)^{1-p} S\bigr)^{1/p p'} {\cal K}^{1/p} \| \, t' \, \|_{L_p(\Bbb{R}_+)}.$$
It remains to apply  (\ref{k28a}). $\square$

\medskip

{\bf Remark 5.} If the requirement
$$
u=0\quad on\; a\; neighbourhood\; of\; a\; closed\; subset\; E\; of\; \overline \Omega
$$
is added in the formulation of Theorems 4 and 5, the same proofs give conditions for the validity of (\ref{k41}),
similar to (\ref{k42}) and (\ref{k45}). The only new feature is the a restriction
$$
\Omega \cap \partial (\Omega \backslash F_2)\quad
is\; at\; a\; positive\; distance\; from\; E.
$$
In the important particular case $E=\partial \Omega $, which corresponds to zero Dirichlet data on $\partial \Omega $,
the conditions (\ref{k42}) and (\ref{k57}) become
\begin{equation}\label{k64}
\mu (F, \Omega \backslash G)^{p/q}\leq B\; {\rm cap}_p(F; G)
\end{equation}
and
\begin{equation}\label{k65}
\mu (F, \Omega \backslash G)\leq \nu \bigl ( ({\rm cap}_p(F; G) )^{1-p}\bigr ),
\end{equation}
respectively, where $F$ is closed and $G$ is open, $G\supset F$, and the closure of $G$ is compact and situated in $\Omega $. The capacity ${\rm cap}_p(F; G)$ is defined by (\ref{k22}) with  $\Omega = G$.

\smallskip

Using lower estimates for the $p$-capacity in terms of area minimizing functions, one obtains sufficient conditions from (\ref{k42}), (\ref{k45}) (\ref{k64}) and (\ref{k65}) formulated in geometrical terms in the spirit of Corollary 2. For example, by (\ref{k64}) and (\ref{k65}), inequalities (\ref{k42}) and (\ref{k57}) hold for  all $u\in C_0^\infty(\Omega)$ if, respectively
$$\mu (F, \Omega \backslash G)^{p/q}\leq B \Bigl(\int_{m_n(F)}^{m_n(\Omega\backslash G)}\frac{dv}{\lambda(v)^{p/(p-1)}} \Bigr)$$
and
$$\mu (F, \Omega \backslash G) \leq \nu\Bigl(\int_{m_n(F)}^{m_n(\Omega\backslash G)}\frac{dv}{\lambda(v)^{p/(p-1)}} \Bigr),$$
where $F$ and $G$ are the same as in (\ref{k64}) and (\ref{k65}). $\square$

\smallskip

By obvious modifications of the proof of sufficiency in Corollary 4 one deduces  the following assertion from Theorem 5.

\smallskip

{\bf Corollary 6.}  (One-dimensional case) {\it With the notation used in Corollary 5, suppose that
$$\mu(I, \Omega\backslash J) \leq \nu(r).$$
Then there exists a positive constant $c$ depending only on $p$ and such that
$$\langle u\rangle_{p,\mu} \leq c\, S^{1/pp'} K^{1/p}\, \|u'\|_{L_p(\Omega)}$$
for all} $u\in C^\infty(\Omega)$.

{\bf Remark 6.} Let us show that the condition $K<\infty$, 
which appeared in Theorem 5, is sharp. Suppose that there exists a positive constant $C$ independent of $u$ and such that 
\begin{equation}\label{71f}
\int_{\Bbb{R}}\int_{\Bbb{R}} |u(t) - u(\tau)|^p\,  \nu'' (t-\tau) dt \, d\tau \leq C \int_{\Bbb{R}} |u'(t)|^p dt,
\end{equation}
where $\nu$ is a convex function in $C^2(\Bbb{R})$. We take  an arbitrary $N>0$ and put $u(t) =\min \{|t|, N\}$ into (\ref{71f}). Then 
$$\int_0^{N/2} \int_\tau^N (t-\tau)^p \,  \nu'' (t-\tau) dt \, d\tau \leq 2CN$$
and setting here $t=\tau +s$, we obtain
$$\frac{1}{2} pN \int_0^{N/2} s^{p-1}\, |\nu'(s)| \, ds \leq p \int_0^{N/2} \int_0^{N-\tau} s^{p-1} |\nu'(s)| ds\, d\tau \leq 2CN.$$
Hence $K\leq 4 p^{-1} C$.

\medskip

{\bf Remark 7.} It seems appropriate, in conclusion, to say a few words about the lower estimate for the difference seminorm $\langle u\rangle_{p,\mu}$, similar to the classical Sobolev inequality:
\begin{equation}\label{72f}
\Bigl( \int_\Omega |u|^q\, \nu(dx) \Bigr)^{1/q} \leq C \, \langle u\rangle_{p, \mu},
\end{equation}
where $\Omega$ is a subdomain of a Riemannian manifold, $\mu$ and $\nu$ are measures in $\Omega \times \Omega$ and $\Omega$, respectively, and $u$ is an arbitrary function in $C_0^\infty(\Omega)$. Suppose that $q\geq p\geq 1$. Then a condition, necessary and sufficient for (\ref{72f}), is the isocapacitary inequality
\begin{equation}\label{73f}
\sup \limits_{\{F\}} \frac{\nu(F)^{p/q}}{{\rm cap}_{p,\mu} (F; \Omega)} <\infty ,
\end{equation}
where $F$ is an arbitrary compact set in $\Omega$ and 
the capacity is defined by
$${\rm cap}_{p,\mu} (F; \Omega) = \inf \bigl\{ \langle u\rangle_{p,\mu}^p: \, u\in C_0^\infty(\Omega), \, u\geq 1 \,\, {\rm on} \,\, F \bigr\}.$$

\smallskip

The necessity of (\ref{73f}) is obvious and the sufficiency results directly from the inequality
\begin{equation}\label{74f}
\int_0^\infty {\rm cap}_{p,\mu} (N_t; \Omega) \, d(t^p) \leq c(p) \,  \langle u\rangle_{p,\mu}^p
\end{equation}
(see \cite{M12} for the proof and history of (\ref{74f})).

\smallskip

Although providing a universal characterization of (\ref{72f}), the condition (\ref{73f}) does not seem satisfactory when dealing with concrete measures and domains. This is related even to    one-dimensional case (cfr. Problem 2 \cite{KMP} ). As an example  of a more visible criterion, consider the measure $\mu$ on $\Bbb{R}^n\times \Bbb{R}^n$ given by 
\begin{equation}\label{75f}
\mu(dx, dy) = |x-y|^{-n-p\alpha} dx\, dy
\end{equation}
with $0<\alpha <1$ and $\alpha p<n$. This measure generates a  seminorm in the homogeneous Besov space $b_p^\alpha (\Bbb{R}^n)$. With this particular choice of $\mu$, we have by Theorem 8.7.1 and Remark 8.6/3 \cite{M10} that (\ref{72f})  holds with $q>p>1$ and $q\geq p=1$ if and only if 
\begin{equation}\label{76f}
\sup\limits_{x\in \Bbb{R}^n, \, \rho>0}\frac{
\nu(B(x,\rho))^{p/q}} {\rho^{n-p\alpha}} <\infty.
\end{equation}
 The inequality (\ref{76f}) is the same as (\ref{73f}) with $F= B(x,\rho)$. It is unknown whether the replacement of arbitrary sets by balls in (\ref{73f})  is possible for the general $\mu$ and $\Omega = \Bbb{R}^n$ in  (\ref{73f}). If not, what are sharp requirements allowing this replacement?

\smallskip

Let $q=p>1$, $\Omega = \Bbb{R}^n$ and let $\mu$ be given by (\ref{75f}). Then (\ref{72f})  holds simultaneously with the inequality
\begin{equation}\label{77f}
\int_{\Bbb{R}^n} |u|^p\, \nu(dx) \leq c \, \| (-\Delta)^{\alpha/2} u\|_{L_p(\Bbb{R}^n)}^p
\end{equation}
because both (\ref{72f}) and (\ref{77f}) are equivalent to isocapacitary inequalities of the type (\ref{73f})  with equivalent capacities in the right-hand side (see \cite{AH}, Sect.4.4).

\smallskip

Note that (\ref{77f}) is the so called trace inequality for the Riesz potential operator $I_\alpha : = (-\Delta)^{-\alpha/2}$. This inequality  has been studied intensively (see \cite{Ve} for a survey of this area). First of all, the simplest estimate
$$\nu(B) \leq c\, m_n(B)^{1-p\alpha/n} \quad {\rm for}\,\, {\rm all}\,\, {\rm balls}\,\, B,$$
being necessary for (\ref{77f}), is not sufficient for it (see \cite{Ad} and \cite{AH}). However, there exist other conditions involving no capacity, which are necessary and sufficient for  (\ref{77f}). They are as follows:

(i) For every ball $B$,
 $$ \quad \int_B (I_\alpha \nu_B)^p dx \leq c\, \nu(B),$$
where $\nu_B$ be the restriction of $\nu$ on $B$, see \cite{KS}.

(ii) Almost everywhere in $\Bbb{R}^n$,
$$I_\alpha(I_\alpha \nu)^{p'} \leq c \, I_\alpha \nu,$$
see \cite{MV}.

(iii) For every dyadic cube $P$ of side length $ \ell(P)$,
$$\sum_{Q\subset P} \bigl( \nu(Q)\, \ell(Q) ^{\alpha -n/p} \bigr)^{p'} \leq c\, \nu(P),$$
where the sum is taken over all dyadic cubes $Q$ contained in $P$, see \cite{Ve}, Sect. 3.

\smallskip

In accordance with the equivalence of  (\ref{72f}) and (\ref{77f}) mentioned previously, the criteria (i)-(iii) characterize not only (\ref{77f}) but also (\ref{72f}) with $q=p>1$ and $\mu$ defined by (\ref{75f}). It is unclear how these criteria could be modified to characterize (\ref{72f}) with an arbitrary  $\mu$.


\bigskip

\begin{thebibliography}{RRRR}

\bibitem{Ad} D.R. Adams, {\it My love affair with the Sobolev inequality}, Sobolev Spaces in Mathematics  I, Sobolev Type Inequalities, 2008, Springer, p. 1-24.

\bibitem{AH} D.R. Adams, L.-I. Hedberg, {\it
Function Spaces and Potential Theory}, Springer, 1996.


\bibitem{ACMM} A. Alvino, A. Cianchi, V. Maz'ya, A. Mercaldo, {\it Well-posed elliptic Neumann problems involving irregular data and domains}, to appear.

\bibitem{Bl} G.A. Bliss,  {\it An integral inequality},  J. London Math. Soc. {\bf 5} (1930), 40--46.


\bibitem{BM} Yu. Burago,  V.  Maz'ya,  {\it Certain Questions of
Potential
Theory and Function Theory for Regions
with Irregular Boundaries}. (Russian) Zap. Nau\v cn. Sem. Leningrad.
Otdel.
Mat. Inst. Steklov. (LOMI) {\bf 3} 1967;
English translation:  Potential Theory and Function Theory for
Irregular
Regions.
 Seminars in Mathematics, V. A. Steklov Mathematical
Institute,
Leningrad, Vol. 3 Consultants Bureau, New York 1969.


\bibitem{Ch} J. Cheeger, {\it A lower bound for the smallest eigenvalue of the Laplacian}, Problems in Analysis, R. Gunning ed., Princeton U.P., 1970, 195--199.


\bibitem{CM} A. Cianchi,  V. Maz'ya,  {\it Neumann problems and isocapacitary inequalities},  J. Math. Pures Appl. {\bf (9) 89} (2008), no. 1, 71--105.

\bibitem{Gr} A. Grigor'yan,  {\it Isoperimetric inequalities and capacities on Riemannian manifolds}, In the book: Operator Theory, Advances and Applications, Vol. 110, The Maz'ya Anniversary Collection, vol. 1, Birkh\"auser, 1999, 139--153.

\bibitem{KS}   R. Kerman, E.  Sawyer, 
 {\it  The trace inequality and eigenvalue estimates 
for Schr\"oding\-er operators}, 
 Ann. Inst. Fourier (Grenoble)  {\bf 36}  (1986), 
 207-228. 


\bibitem{KK} J. Kinnunen, R. Korte, {\it Characterization of Sobolev inequalities on metric spaces}, arXiv: 0709.2013v1 [mathAP].

\bibitem{KMP} A. Kufner, L. Maligranda, L.-E. Persson, {\it The Hardy Inequality, About its History and Some Related Results}, Pilsen, 2007.

\bibitem{M1} V.  Maz'ya,  {\it Classes of domains and imbedding theorems for
function spaces}. Sov. Math. Dokl. {\bf 1} (1960),
882--885.

\bibitem{M2} V.  Maz'ya, {\it On $p$-conductivity and theorems on embedding
certain functional spaces into a $C$-space}, Sov. Math. 
 Dokl.   {\bf 2} (1961), 1200--1203.

\bibitem{M2a} V. Maz'ya,  {\it Some estimates of solutions of second-order
elliptic equations},  Sov. Math.  Dokl. 
 {\bf 2} (1961), 413--415.


\bibitem{M3} V.  Maz'ya, {\it The negative spectrum of the
higher-dimensional
Schr\"odinger operator},  Sov. Math. 
 Dokl.   {\bf 3} (1962), 808--810.

\bibitem{M4} V.  Maz'ya, {\it On the solvability of the Neumann problem}.
  Sov. Math.  Dokl.  {\bf 3} (1962),
1595--1598.

\bibitem{M5} V.  Maz'ya, {\it The Dirichlet problem for elliptic equations
of
arbitrary order in unbounded domains}, Sov. Math. 
 Dokl.  {\bf 4} (1963), 860--863.


\bibitem{M6} V. Maz'ya, {\it On the theory of the higher-dimensional
Schr\"odinger operator}, 
 Izv. Akad. Nauk SSSR Ser. Mat. {\bf 28}, no. 4 (1964), 1145--1172.


\bibitem{M7} V.  Maz'ya, {\it The Neumann problem in regions with
nonregular
boundaries},   Sibirsk.
Mat. \v Z. {\bf 9} (1968), 1322--1350. 

\bibitem{M8} V.  Maz'ya, {\it Weak solutions of the Dirichlet and Neumann
problems}. Trudy Moskov. Mat. Ob\v s\v c. {\bf 20}  (1969),
137--172.

\bibitem{M9} V.  Maz'ya, {\it Certain integral inequalities for functions
of
several variables}, 
 Problems of mathematical analysis, no. 3: Integral and
differential operators, Differential equations (Russian), pp.
33--68. Izdat. Leningrad. Univ., Leningrad, 1972; English translation: J. Soviet Math. {\bf 1} (1973), 205--234.

\bibitem{M10} V.  Maz'ya,   {\it Sobolev Spaces}, Springer, 1985.

\bibitem{M11} V.  Maz'ya, {\it Lectures on isoperimetric and isocapacitary inequalities 
in the theory of Sobolev spaces}, 
{Contemporary Mathematics}, 
{\bf 338}, {Heat Kernels and Analysis on Manifolds, Graphs, and Metric
Spaces}, 
{American Math. Society}, 2003, 
{307--340}.


\bibitem{M12} V.  Maz'ya, 
{\it Conductor and capacitary inequalities for functions on
topological
spaces and their applications to Sobolev type imbeddings}, 
  {Journal of Functional Analysis}, {224}, no. {2} (2005), 408--430.

\bibitem{M13} V.  Maz'ya, {\it Conductor inequalities and criteria for Sobolev type two-weight imbeddings}, 
{J. Comput. Appl. Math.}, 
{\bf 194}, no. {1} (2006), 
{94--114}.


\bibitem{M14} V.  Maz'ya, 
{\it Analytic criteria in the qualitative spectral analysis of the Schr\"odinger operator}, 
{ Spectral Theory and Mathematical Physics: a Festschrift in honor of Barry Simon's 60th birthday,   Proc. Sympos. Pure Math.,  Part 1, Amer. Math. Soc.,  Providence, RI}
{\bf 76}, 2007, 
{257--288}.


\bibitem{MV} V. Maz'ya,  I.   Verbitsky, 
 {\it Capacitary
estimates for
fractional integrals, with applications to partial differential
equations and Sobolev multipliers},  Arkiv f\"or Matem. {\bf 33}
 (1995),   81-115.


\bibitem{PS} G. P\'olya, G. Szeg\"o, {\it Isoperimetric Inequalities in Mathematical Physics},  Ann. Math. Stud. 27, Princeton Univ. Press, Princeton, 1951.


\bibitem{St} E.M. Stein, {\it Singular Integrals and Differentiability Properties of Functions}, Princeton University Press, Princeton, N. J., 1970.

\bibitem{Ve} I. Verbitsky, {\it  Nonlinear potentials and trace inequalities},  Operator Theory, Advances and Applications, Vol. 110, The Maz'ya Anniversary Collection, vol. 2, Birkh\"auser, 1999, 323 - 343.



\end{thebibliography}
\end{document}